\documentclass{article}

\usepackage[utf8]{inputenc}

\usepackage{dsfont}

\begin{document}

\noindent {\large \bf \textsf{A Minimalist Proof of the Law of Quadratic Reciprocity}} \\

\noindent This proof is similar to one of the proofs of Eisenstein [1], albeit shorter.\\

\noindent Let $S_n(t)=\sum_{x_1+x_2+\dots+x_n=t} \left(\frac{x_1x_2\cdots x_n}{p}\right)$, where $t, x_i \in \mathbf{Z}_p$ and $\left( \frac{a}{p}\right)$ is the Legendre symbol. Note that for an odd $n$ and any $a \in \mathbf{Z}_p/\{0\}$, 
\[   S_n(t)=\left(\frac{a}{p}\right)^n\sum_{\frac{x_1}{a}+\cdots+\frac{x_n}{a}
=\frac{t}{a}} \left(\frac{\frac{x_1}{a}\cdots\frac{x_n}{a}}{p}\right)=\left(\frac{a}{p}\right)S_n(t/a).  \] 
In particular, substituting $t=0$ and $a$ a nonresidue, we get $S_n(0)=0$, and setting $t=a$, we get $S_n(a)=(\frac{a}{p})S_n(1)$; and by similar reasoning, $S_2(a)=(\frac{a}{p})^2S_2(1)=S_2(1)$. Furthermore, $S_2(0)=\sum_{x \in \mathbf{Z}_p/\{0\}} (\frac{-x^2}{p})=(\frac{-1}{p})(p-1)$, and, recalling that $\sum_{t \in \mathbf{Z}_p}(\frac{t}{p})=0$,  \[S_2(1)=\sum_{x \in \mathbf{Z}_p/\{0\}} \left(\frac{x^2x^{-1}(1-x)}{p}\right)=\sum_{x \in \mathbf{Z}_p/\{0\},\ y=x^{-1}} \left(\frac{y-1}{p}\right)=-\left(\frac{-1}{p}\right).\] Therefore, for odd $n$ we have $S_{n+2}(1)=\sum_{t \in \mathbf{Z}_p/\{0, 1\}} S_n(t)S_2(1-t) + S_n(1)S_2(0) + S_n(0)S_2(1) = \sum_{t \in \mathbf{Z}_p/\{0, 1\}} (\frac{t}{p}) S_n(1) (-(\frac{-1}{p})) + S_n(1)(\frac{-1}{p})(p-1) = S_n(1)(\frac{-1}{p})p.$ Since $S_1(1)=(\frac{1}{p})=1$, by induction for odd $n$ $S_n(1)=((\frac{-1}{p})p)^{\frac{n-1}{2}}=p^{\frac{n-1}{2}}(-1)^{\frac{p-1}{2}\frac{n-1}{2}}$. Therefore, by Euler's criterion, for an odd prime $q$, $S_q(1) \equiv (\frac{p}{q})(-1)^{\frac{p-1}{2}\frac{q-1}{2}} \pmod q$. Now, the $q$ cyclic shifts of a given $q$-tuple $(x_1, \dots, x_q)$ are distinct unless all $x_i$ are equal, since the period of its repeated single-position cyclic shift divides $q$, and so is $q$ or 1. When they are distinct, their total contribution to the sum defining $S_q(1)$ is $q(\frac{x_1x_2\cdots x_q}{p})$, which is divisible by $q$. Therefore, modulo $q$ (we take $q\neq p$), \[S_q(1) \equiv \sum_{x_1+x_2+\dots+x_q=1,\ x_1=x_2=\dots=x_q} \left(\frac{x_1x_2\cdots x_q}{p}\right)=\left(\frac{q^{-1}}{p}\right)^q=\left(\frac{q}{p}\right).\] So $(\frac{p}{q})(-1)^{\frac{p-1}{2}\frac{q-1}{2}}$ and $(\frac{q}{p})$ are congruent to $S_q(1)$, and thus to each other, modulo $q$ --- but they both are numbers of the form $\pm 1$, so they are equal, which is the law of quadratic reciprocity.

\rightline{---Submitted by Bogdan Veklych, University of Wisconsin--Madison}

\bigskip
\footnoterule
\footnotesize{doi.org/10.XXXX/amer.math.monthly.122.XX.XXX}


\end{document}